\documentclass{amsart}

\usepackage[dvips]{graphics}
\usepackage{amsthm,amssymb,latexsym,abbreviations, indentfirst}
\usepackage{epsfig}
\usepackage{color}
\usepackage[mathscr]{eucal}

\newtheorem*{theorem}{Theorem}

\newtheorem{corollary}{Corollary}
\newtheorem{lemma}{Lemma}
\newtheorem{fact}{Fact}
\newtheorem*{claim}{Claim}

\theoremstyle{definition}
\newtheorem{definition}{Definition}

\theoremstyle{remark}
\newtheorem{remark}{Remark}

\newcommand{\PPP}{\mathbb{P}}
\newcommand{\EEE}{\mathbb{E}}

\newcommand{\eh}{excess height }

\newcommand{\nid}{\noindent}

\renewcommand{\qedsymbol}{\ensuremath{\blacksquare}}

\begin{document}

\title[Diffusion in random environment]
{Diffusion in random environment and \\ the renewal theorem}
\author{Dimitrios Cheliotis}
\address{Department of Mathematics\\
University of Toronto\\100 St. George St\\
Toronto, ON M5S 3G3 \\ Canada} \email{dimitris@math.toronto.edu}
\urladdr{http://www.math.toronto.edu/dimitris/}
\thanks{Research partially supported by an anonymous Stanford Graduate Fellowship, A.Onassis Scholarship, and by NSF grants
\#DMS-0072331, \#DMS-FRG-0244323} 
\date{August 8, 2004}

\begin{abstract}
According to a theorem of S. Schumacher and T. Brox, for a
diffusion $X$ in a Brownian environment it holds that
$(X_t-b_{\log t})/\log^2t\to 0 $ in probability, as $t\to\infty$,
where $b$ is a stochastic process having an explicit
description and depending only on the environment. We compute the distribution of the number of sign changes
for $b$ on an interval $[1,x]$ and study some of the consequences
of the computation; in particular we get the probability of $b$
keeping the same sign on that interval. These results have been
announced in 1999 in a non-rigorous paper by P. Le Doussal, C.
Monthus, and D. Fisher and were treated with a Renormalization
Group analysis. We prove that this analysis can be made rigorous
using a path decomposition for the Brownian environment and
renewal theory. Finally, we comment on the information these results give about the behavior of the diffusion.
\end{abstract}

\maketitle

\section{Introduction}\label{S:intro}

On the space $\C{W}:=C(\D{R})$ consider the topology of
uniform convergence on compact sets, the corresponding $\gs$-field
of the Borel sets, and $\PPP$ the measure on $\C{W}$ under which
the coordinate the processes $\{w(t):t\ge0\},\{w(-t):t\ge0\}$ are
independent standard Brownian motions.

Also let $\Omega:=C([0,+\infty))$, and equip it with the
$\gs$-field of Borel sets derived from the topology of uniform
convergence on compact sets. For $w\in\C{W}$, we denote by
$\textbf{P}_w$ the probability measure on $\Omega$ such that
$\{\omega(t):t\ge0\}$ is a diffusion with $\omega(0)=0$ and
generator
\[
\frac{1}{2}e^{w(x)}\frac{d}{dx}\left(e^{-w(x)}\frac{d}{dx}\right).
\]
The construction of such a diffusion is done with scale and time
transformation from a one-dimensional Brownian motion (see, e.g.,
\cite{SC1}, \cite{SH1}). Using this construction, it is easy to
see that for $\PPP$-almost all $w\in\C{W}$ the diffusion does not
explode in finite time; and on the same set of $w$'s it satisfies
the formal SDE

\begin{equation}\label{FSDE}
\begin{array}{rl}
d\omega(t)=&d\beta(t)-\frac{1}{2}w'(\omega(t))dt, \\
\omega(0)=&0,
\end{array}
\end{equation}
where $\beta$ is a one-dimensional standard Brownian motion.

Then consider the space $\C{W}\times\Omega$, equip it with the
product $\gs$-field, and take the probability measure defined by
\[
d\C{P}(w,\omega)=d\textbf{P}_w(\omega)d\PPP(w).
\]
The marginal of $\C{P}$ in $\Omega$ gives a process that is known
as diffusion in a random environment; the environment being the
function $w$.

S. Schumacher (\cite{SC1, SC2}) proved the following result.

\begin{fact} \label{SchumProp}
There is a process $b:[0,\infty)\times\C{W}\to\D{R}$ such that
for the formal solution $\omega$ of \eqref{FSDE} it holds
\begin{equation}\label{Schum}
\frac{\omega_t}{(\log t)^2}-b_1(w^{(\log t)})\to 0 \text{ in
$\C{P}$ as $t\to+\infty$,}
\end{equation}
where for $r>0$ we let $w^{(r)}(s)=r^{-1}w(sr^2) \text{ for all
$s\in\D{R}$} .$
\end{fact}
We will define the process $b$ soon.
 This result shows the dominant
effect of the environment, through the process $b$, on the asymptotic behavior of the
diffusion. The results we prove in this paper concern
the process $b$. In subsection \ref{connection} we commend on their implications
for the behavior of the diffusion itself.

Besides this diffusion model, there is a discrete time and space
analog, known as Sinai's walk, which was studied first. Sinai's
pioneering paper \cite{SI} identified the role of the process $b$
in the analogous to \eqref{Schum} limit theorem for the walk. Then
S. Schumacher proved in \cite{SC2} (see also \cite{SC1} for the
results without the proofs) a more general statement than the
above proposition where the environment $w$ was not necessarily a
two sided Brownian motion,  while T. Brox \cite{BR} gave a
different proof in the Brownian case. H. Kesten \cite{KE} computed
the density of $b_1$ in the case we consider, and Tanaka
\cite{TA1} generalized the computation to the case that $w$ is a
two sided symmetric stable process. Localization results have been
given for the Sinai walk by Golosov (\cite{GO}, actually, for the
reflected walk) and for the diffusion model by Tanaka \cite{TA2}.
Also Tanaka (\cite{TA3}, \cite{TA1}) studied the cases where the
environment is non-positive reflecting Brownian motion,
non-negative reflecting Brownian motion, or Brownian motion with
drift. Finer results on the asymptotics of Sinai's have been
obtained
 by Z. Shi and Y. Hu. A survey of some of them as well as a
connection between Sinai's walk and diffusion in random
environment is given in \cite{SH1}. Another connection is
established in \cite{SEI}.

In \cite{DMF}, P. Le Dousal, C. Monthus, and D. Fisher proposed
 a new method for tackling questions related to
asymptotic properties of Sinai's walk, and using it they gave a
host of results.  The method is a Renormalization Group analysis
and it has consequences agreeing with rigorously proved results
(e.g., \cite{DGZ}, \cite{KE}) . This is the starting point of the
present paper. In the context of diffusion in random environment,
we show how one can justify the Renormalization Group method using two tools. The first
is a path decomposition for a two sided standard Brownian motion; the second is the renewal theorem.
Our main results illustrate the use of the method and the way we justify it.

The structure of the paper is as follows. In the remaining of
the introduction we state our results. In Section
\ref{preparation} we provide all the necessary machinery for the
proofs, which are given in Section \ref{proofs}. Some technical
lemmata that we use are proved in Section \ref{lemmata}.

\medskip

\nid We begin by defining the process $b$. \vspace{3pt}

For a
function $w:\mathbb{R}\to\D{R},\ x>0$ and $y_{0}\in\D{R}$ we say
that \textbf{$w$ admits an $x$-minimum at $y_{0}$} if there are
$\ga,\gb\in\mathbb{R}$ with $\ga<y_{0}<\gb$,
$w(y_{0})=\inf\{w(y):y\in[\ga,\gb]\}$ and $w(\ga)\ge w(y_{0})+x$,
$w(\gb)\ge w(y_{0})+x$. We say that \textbf{$w$ admits an
$x$-maximum at $y_{0}$} if $-w$ admits an $x$-minimum at $y_{0}$.

\begin{figure}[h]
\begin{center}
\resizebox{6cm}{!}{
\input{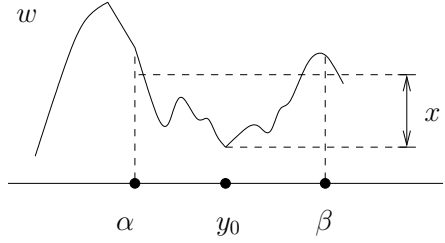}}
\end{center}
\caption{\label{xminP}$w$ admits an $x$-minimum at $y_0$.}
\end{figure}

For convenience, we will call a point where $w$ admits an
$x$-maximum or $x$-minimum an $x$-maximum or an $x$-minimum
respectively.

We denote by $R_x(w)$ the set of $x$-extrema of $w$ and define
\[
\C{W}_1:=\left\{w\in\C{W}: \begin{array}{c}\text{ For every $x>0$
the set $R_x(w)$ has no accumulation point in $\D{R}$,}  \\
\text{it is unbounded above and below,}\\
\text{and the points of $x$-maxima and $x$-minima
alternate.}\end{array}\right\}.
\]
Thus, for $w\in\C{W}_1$ and $x>0$ we can write
$R_x(w)=\{x_k(w,x):k\in \D{Z}\}$ with $(x_k(w,x))_{k\in\D{Z}}$
strictly increasing, $x_0(w,x)\le0<x_1(w,x)$ ,
$\lim_{k\to-\infty}x_k(w,x)=-\infty, \quad
\lim_{k\to\infty}x_k(w,x)=\infty$. It holds that
$\PPP(\C{W}_1)=1$, and the easy proof of this fact is given in
Lemma \ref{discrete}.

\begin{definition}

The process $b:[0,+\infty)\times\C{W}\to \D{R}$ is defined for
$x>0$ and $w\in \C{W}_1$ as
$$
b_{x}(w):=\begin{cases} x_0(w,x) &  \text{ if $x_0(w,x)$ is an $x$-minimum,} \\
                      x_1(w,x) &  \text{ else,}
        \end{cases}
$$
and $b_x(w)=0$ if $x=0$ or $w\in\C{W}\setminus\C{W}_1$.
\end{definition}

\begin{remark}
 In the definition of $b_x(w)$ we do
not make use of the entire sequence of $x$-extrema. The reason we
introduce this sequence is that we plan to study the evolution
of the process $b$ as $x$ increases. Since $R_{\tilde x}(w)\subset
R_x(w)$ for $x<\tilde x$, the later values of $b_\cdot(w)$ are
elements of $R_x(w)$. For $\tilde x$ large enough, the points
$x_0(w,x)$, $x_1(w,x)$ will not be $\tilde x$-extrema.
\end{remark}

\begin{remark}
We will decompose the process $w$ at the endpoints of the
intervals $\{[x_k(w,x), x_{k+1}(w,x)]:k\in\D{Z}\}$ and study its
restriction to each of them.  Of course $[x_0(w,x), x_1(w,x)]$ has
a particular
 importance for the process $b$; and it is in the study of $w|[x_0(w,x), x_1(w,x)]$ that the renewal theorem enters
 (see Lemma \ref{cslope}).
\end{remark}

\begin{remark}\label{bscaling}
It is clear that $b$ satisfies $b_x(\frac{1}{a}w(c \cdot))=\frac{1}{c}b_{ax}(w)$ for all $a, c, x>0$, and $w\in\C{W}$.
So that the quantity $b_1(w^{(\log t)})$ appearing in \eqref{Schum} equals also $b_{\log t}(w)/(\log t)^2$.
\end{remark}

\subsection{Sign changes of $b$}

For $x\ge1$ define on $\C{W}_1$ the random variable
$$ k(x)=\#\text{~times~~} b_\cdot(w) \text{~~has changed sign in~}
[1,x].$$ The main result of the paper is the computation of the
generating function of $k(x)$.

\begin{theorem}\label{thm1}
For $x\ge1$ and $z\in\D{C}$ with $|z|<1$, it holds
\begin{equation}\label{genfun}
\EEE(z^{k(x)})=c_1(z)\ x^{\gl_1(z)}+c_2(z)\ x^{\gl_2(z)},
\end{equation}
where $$\gl_{1}(z)=\frac{-3+\sqrt{5+4z}}{2},\quad
\gl_{2}(z)=\frac{-3-\sqrt{5+4z}}{2},$$ and
$$c_{1}(z)=((z-1)/3-\gl_2(z))/(\gl_{1}(z)-\gl_{2}(z)), \quad
 c_{2}(z)=(-(z-1)/3+\gl_1(z))/(\gl_{1}(z)-\gl_{2}(z)).$$
\end{theorem}
From this we extract several corollaries. Corollary \ref{nozero} and Corollary \ref{nzeros} are immediate,
while the rest require some work and are proved in section \ref{proofs}.

\medskip

\nid Corollary \ref{nozero} follows by taking $z\to 0$ in \eqref{genfun}.
\begin{corollary} \label{nozero}
\begin{equation}
 \PPP(b_{\cdot}(w) \text{~~~keeps the same sign in~~} [1,
 x])/
 x^{(-3+\sqrt{5})/2}\to 1/2+7\sqrt{5}/30
\end{equation}
as $x\to +\infty$.
\end{corollary}

\begin{corollary}\label{Markov}
The increasing process of points $(X_k)_{k\ge1}$ where $b$ changes sign in $[1, +\infty)$
has the form $X_k=X_1r_1\cdots r_{k-1}$, $k\ge2$, where $X_1$ is the smallest such point
and the $r_i$'s are i.i.d with density
\begin{equation} \label{rdensity}
 f(r)=\frac{1}{\gl_1-\gl_2}(r^{\gl_1-1}-r^{\gl_2-1}), \quad
r\ge1,
\end{equation}
where $\gl_1=\gl_1(0), \gl_2=\gl_2(0)$.
\end{corollary}
\nid Now observe that
$$k(t)=\sup\{n\in\D{N}:X_n\le t\}=\sup\{n\in\D{N}:\log X_1+\log r_1+\cdots +\log r_{n-1} \le \log t\}.$$
Since $\EEE(\log r_1)=3$ and $\log X_1$ is finite a.s. (e.g., by
Corollary \ref{nozero}, $\EEE(\log X_1)<+\infty$), the next statement follows from renewal theory.

\begin{corollary} \label{nzeros} $k(t)/\log t\rightarrow 1/3$ as
$t\rightarrow +\infty$  $\PPP$-a.s.
\end{corollary}
Corollary \ref{Markov} allows the following strengthening of the above theorem.
\begin{corollary} \label{strenghening}
Relation \eqref{genfun} holds for all $z\in\D{C}\setminus(-\infty,-5/4].$
\end{corollary}
For $t\in(0, \infty)$ consider the random variable $k(e^t)/t$
and let $\mu_t$ be its distribution measure. Then the following
holds.
\begin{corollary}\label{LDP1} The family of measures $(\mu_t)_{t>0}$
satisfies a Large Deviation Principle with speed $t$ and good
rate function
$$ I(x)=
\left\{
\begin{array}{lr} x~\log\big(2x(x+\sqrt{x^2+5/4})\big)+\frac{3}{2}-(x+\sqrt{x^2+5/4})
 &\quad  \text{ if }x\in[0,+\infty),  \\
+\infty  &\quad   \text{ if } x\in(-\infty,0). \end{array}
    \right.
  $$
\end{corollary}
In \cite{DMF}, Corollary \ref{Markov} appears in paragraph
IV.B with a different justification. We state it here because we
need it for the proof of Corollary \ref{LDP1}. The large deviation
result of Corollary \ref{LDP1} is the precise mathematical
interpretation of the discussion in paragraph IV.A of the same paper.

In Corollary \ref{nozero}, the condition $k(x)=0$ means that the
process $b$ tends to keep the diffusion away form zero in the time
interval $[e,e^x]$ (since the diffusion localizes around $b$, and
$b$ keeps sign on $[1,x]$). For the event that the diffusion
\textit{hits} zero there are two interesting relevant papers. The
first one is by Y. Hu \cite{HU} who treats the annealed
asymptotics of the first time of hitting zero after time $t$ as
$t\to+\infty$ . The second is by F. Comets and S. Popov \cite{CP}
and refers to a related model. That is, it considers a process
$(X_t)_{t\ge0}$ on $\D{Z}$ that runs in continuous time in an
environment $\omega$ satisfying the conditions of the Sinai model
and studies, among other things, the asymptotics of the quenched
probability $\textbf{P}_w(X_t=0|X_0=0)$  as $t\to+\infty$.

\subsection{The process $b$ and the diffusion} \label{connection}
The results of the previous subsection concern the process $b$, which is a functional of the environment.
And our motivation for studying $b$ was the localization results involving this process
(the simplest being Fact \ref{SchumProp}, and  keep in mind Remark \ref{bscaling}).
An obvious question is what we can infer about the behavior of the diffusion from our results.

Using the representation of the diffusion as a time and scale change of Brownian
motion, one can show easily that the diffusion is recurrent and 0 is a regular point for
$(0, +\infty)$ and $(-\infty, 0)$.
Consequently, for all $c>0$, the diffusion visits 0 in the time interval $[c, +\infty)$
and in its first visit there it scores an infinite number of sign changes.
So there can be no direct connection with the corresponding number for the process $b$.
One can consider, say, the number of sign changes of the diffusion between times where the diffusion
achieves a positive or negative record value. This number is finite on compact intervals of $(0,+\infty)$
but still not related with the sign changes of $b$ (and it is not hard to see this).

When $b_{\log \cdot}$ jumps to a new value, what happens
is not that just the diffusion goes through that value
shortly before or after that. As it is known (see, e.g., \cite{BR}, \cite{TA5}), the impact of the jump is that
the diffusion goes to that value and it is trapped in its
neighborhood for a large amount of time. The way to detect the approximate location of such values when we
observe the diffusion in real time (i.e., at time $t$ we know $\omega|[0, t]$)
is to find the site the diffusion has spent the most time thus far.

To make the last statement precise, we use the local time process $\{L_\omega(t,x):t\ge0, x\in\D{R}\}$
that corresponds to the diffusion $\omega$. This process is jointly continuous, and with probability one satisfies
$$\int_0^t
f(\omega(s))ds=\int_{\D{R}}f(x)L_\omega(t,x)dx$$ for all $t\ge0$
and any bounded Borel function $f\in\D{R}^\D{R}$.

For a fixed $t>0$, the set $\mathfrak{F}(t):=\{x\in\D{R}:
L_\omega(t,x)=\sup_{y\in\D{R}}L_\omega (t,y)\}$ of the points with
the most local time at time $t$ is nonempty and compact. Any point
there is called a \textbf{favorite point} of the diffusion at time
$t$. Define $F:(0,+\infty)\to \D{R}$ with
$F(t)=\inf\mathfrak{F}(t)$, the smallest favorite point at time
$t$. Also, for $x>1$ define the interval $I(x):=(x-(\log x)^5,
x+(\log x)^5)$.

In a work in progress we expect to prove the following.
\begin{claim} \label{fplocal}
With $\C{P}$ probability one, there is a strictly increasing
sequence $(t_n(\omega, w))_{n\ge1}$ converging to infinity so that
if we denote by $(x_n(w,t_1))_{n\ge1}$ the sequence of consecutive
values that $b_{\log \cdot}$ takes on the interval $(t_1,
+\infty)$ (remember that $b$ is a step function in any interval
$[x, +\infty)$ with $x>0$), then
$$F((t_n, t_{n+1}))\subset I(x_n) \text{ for all } n\ge 1,$$
and $x_n=b_{\log t}$ for some $t\in (t_n, t_{n+1})$. We
abbreviated $t_n(\omega, w), x_n(w,t_1)$ to $t_n, x_n$.
\end{claim}
Observe that for big $x$, the interval $I(x)$ is a
relatively small neighborhood of $x$. Thus, the claim says that, after some point,
the function $F$ ``almost tracks'' the values of
the process $b_{\log \cdot}$ with the same order and at about the
same time. Consequently, the number
of sign changes of $b_{\log\cdot}$ on an interval $(s,t)$ (with $s,t$ large) would correspond
to the number of sign changes of $F$ on approximately the same
interval. Or, more precisely, the number of sign changes of $F$ on $(s,t)$ and the
corresponding number for $b_{\log\cdot}$ differ by at most two.
It is easy to see that the last statement follows from the Claim above.

\section{Preliminaries} \label{preparation}

As a first step towards the study of the process $b$, we look
at the law of the Brownian path between two consecutive
$x$-extrema as well as the way these pieces are put together to
constitute the entire path.

The first piece of information is provided by Proposition of \S 1
in \cite{NP}.

\begin{fact} \label{renewal}
For every $x>0$, the times of x-extrema of a Brownian motion
$(w_t:t\in\mathbb{R}),\  w_0=0$ build a stationary renewal process
$R_x(w)=\{x_k:k\in\mathbb{Z}\}$ with $(x_k)_{k\in\mathbb{Z}}$
strictly increasing and $x_0\le0<x_1.$ The trajectories between
consecutive x-extrema $(w_{x_k+t}-w_{x_k}: t\in[0, x_{k+1}-x_k]),
k\in\mathbb{Z}$ are independent and for $k\ne0$ identically
distributed (up to changes of sign).
\end{fact}

\begin{figure}[h] \label{P2}
\begin{center}
\resizebox{6cm}{!}{
\input{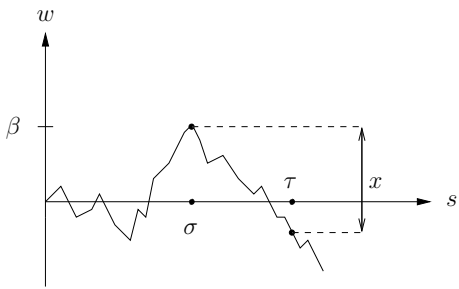}}
\end{center} \caption{The graph of $w$ until $M-w$ hits $x$.}
\end{figure}

In Lemma of \S 1 of the same paper \cite{NP} a description of each
such trajectory is given which we quote (see Figure \ref{P2}).

\nid For $x, t\ge0$ let
\begin{align*}
M_t:&=\sup\{w_s:s\in[0,t]\},\\
\tau:&=\min\{t:M_t=w_t+x\}, \\
\beta:&=M_\tau, \\
\sigma:&=\max\{s\in[0,\tau]:w_s=\beta\}.
\end{align*}

\begin{fact} \label{NPLaplace}
The following hold:
\begin{enumerate}
\item[(i)] The two trajectories $(w_t:t\in[0,\sigma])\text{~~
and~~} (\beta-w_{\sigma+t}:t\in[0,\tau-\sigma])$ are independent.
\item[(ii)]$\beta$ has exponential distribution with mean $x$.
\item[(iii)] The Laplace transform of the law of $\sigma$ given
$\beta$ is
$$\EEE[e^{-s\sigma}|\beta=y]=\exp\{-(y/x)\phi(sx^2)\} \text{\qquad,~s}>0.$$
\item[(iv)] The Laplace transform of the law of  $\tau-\sigma$ is
$$\EEE[e^{-s(\tau-\sigma)}]=\psi(sx^2) \text{\qquad,~s}>0,$$
\end{enumerate}
where $\phi(s)=\sqrt{2s}\coth\sqrt{2s}-1\text{ and
}\psi(s)=\sqrt{2s}/\sinh\sqrt{2s}.$
\end{fact}

\medskip

We call the translation $(w-w(x_k))|[x_k,x_{k+1}]$ of the
trajectory of $w$ between two consecutive $x$-extrema an
\textbf{$x$-slope} (or a slope, when the value of $x$ is clear or
irrelevant), a slope that takes only non-negative values an
\textbf{upward slope}, and a slope taking only non-positive
values a \textbf{downward slope}. We call
$(w-w(x_0))|[x_0,x_1]$ the \textbf{central $x$-slope}.

In the following we will use the operation of ``gluing together''
functions defined on compact intervals. For two functions $f:[\ga,
\gb]\to\D{R}$, $g:[\gamma, \gd]\to\D{R}$, by gluing $g$ to the
right of $f$ we mean that we define a new function $j:[\ga,
\gb+\gd-\gamma]\to\D{R}$ with
\begin{equation*}
j(x)=
 \begin{cases}
f(x) &\mbox{ for } x\in [\ga, \gb], \\
f(\gb)+g(x-\gb+\gamma)-g(\gamma) &\mbox{ for } x\in [\gb,
\gb+\gd-\gamma].
 \end{cases}
\end{equation*}

It is clear that for all $k\neq0$, if $(w-w(x_k))|[x_k,x_{k+1}]$
is an upward $x$-slope, then it is obtained by gluing together a
trajectory of type $(w_t:t\in[0,\sigma])$ to the right
 of a trajectory of type $(\beta-w_{\sigma+t}:t\in[0,\tau-\sigma])$ and then translating
 the resulting path so that it starts at $(x_k, 0)$ .
Similarly for a downward slope.

For any $x$-slope $T:[\alpha,\beta]\to\D{R}$ we call
$l(T):=\beta-\alpha$ , $h(T):=|T(\beta)-T(\alpha)|$ the
\textbf{length} and the \textbf{height} of the slope respectively,
and we denote by $\eta(T):=h(T)-x$ the ``\textbf{excess height}''
of $T$. Also we denote by $|T|$,\ $\theta(T)$, \label{theta} the
slopes with domains $[\alpha,\beta]$, $[0,\beta-\alpha]$ and
values $|T|(\cdot)=|T(\cdot)|$, $\theta(T)(\cdot)=T(\alpha+\cdot)$
respectively.

For any slope $T$, the slope $|\theta(T)|$ is in the set $\C{S}$
defined by

\begin{equation*}
 \C{S}:=
 \left\{f \subset[0,+\infty)^2:\begin{array}{c} f \text{ is a function  such that there is a $l(f)\ge 0$
 with}
 \\ \text{Domain$(f)=[0,l(f)]$, $f$ continuous,}
   \\ \text{ and } 0=f(0)\le f(x) \le f\big(l(f)\big)
\quad \forall x \in [0,l(f)]
\end{array}
\right\}.
\end{equation*}
On $\C{S}$ we define a topology for which the base of
neighborhoods of an element $f\in\C{S}$ is the collection of all
sets of the form $$\{g\in\C{S}:|l(g)-l(f)|<\gep \text{ and $|f(t\
l(f))-g(t\ l(g))|<\gep$ for all $t\in [0,1]$}\}.$$
With this
topology, $\C{S}$ is a Polish space. Equip $\C{S}$ with the Borel
$\gs$-algebra, and define the measures $m_x^r$, $m_x^c$ the first
to be the distribution of
$\theta(|\big(w-w(x_1)\big)|\big|[x_1,x_2])$ and the second to be
the distribution of $\theta(|\big(w-w(x_0)\big)|\big|[x_0,x_1])$
(the superscripts $r$ and $c$ standing for renewal and central).

In the remaining part of this section we compute
the distribution of the length and excess height of a slope picked from $m_x^r$ or $m_x^c$.
We assume $x=1$ since the scaling property of
Brownian motion gives the corresponding results for the case
$x\ne1$.

Let $T$ be a slope picked from $m_1^r$. Earlier we described the
way an upward slope is formed. Its length is the sum of two
independent independent random variables $Z_1, Z_2$ with
$Z_1\overset{law}{=}\gs, Z_2\overset{law}{=}\tau-\gs$ (where we
take $x=1$ in their definitions just before Fact \ref{NPLaplace}).
But from (i) of Fact \ref{NPLaplace}, $\gs$ and $\tau-\gs$ are
independent. Thus, $l(T)\overset{law}{=}\tau$. By definition,
$\tau$ is the time where the reflected process $M-w$ hits one.
This reflected process has the same law as $|w|$, and the
Laplace transform of the time it first hits one is known as
\begin{equation} \label{hitting1}
\EEE[e^{-\gl l(T)}]=(\cosh \sqrt{2\lambda})^{-1} \text{ for } \gl>0.
\end{equation}
Also $\EEE\left(l(T)\right)=\EEE(\tau)=1$.
Using the Laplace inversion formula (see \cite{MH}, pg 531) we find
that the density of $l=l(T)$ is
\begin{equation} \label{preKesten}
f_l(x)=\frac{\pi}{2}\sum\limits_{k\in\mathbb{Z}}(-1)^k(k+\frac{1}{2})\exp[-\frac{\pi^2}{2}(k+\frac{1}{2})^2~x~]
\quad ,~x>0.
\end{equation}
We note for future reference that, by (ii) of Fact
\ref{NPLaplace}, for any $a>0$, the excess height of a slope picked
from $m_a^r$ is exponential with mean $a$; i.e., it has density
\begin{equation} \label{expona}
p_a(x)=a^{-1}e^{-x/a} \qquad ,~x>0.
\end{equation}
Now let $T_0$ be a slope picked from $m_1^c$. More specifically, take $T_0$ to be the central 1-slope.
Observe
that, by virtue of Fact \ref{renewal}, we can ``start a renewal process at
$-\infty$'' with i.i.d. alternating upward and downward slopes, and ask what the
characteristics are of the slope covering zero. The renewal
theorem says that the length of the slope covering zero is picked
from the distribution of $l$ given in \eqref{preKesten} with
size-biased sampling. Once the length, say $z$, is picked, we
expect that the remaining characteristics of the slope, ignoring
direction (upward or downward), are determined by the law of
$T|l(T)=z$ under $m_1^r$. We give a formal proof of this. Notice
that a regular conditional distribution for the random variable
$T$ given the $\sigma$-field $\sigma(l(T))$ exists because the
space $\C{S}$ is Polish.

\begin{lemma} \label{cslope}
For any measurable subset A of $\C{S}$, it holds
\begin{equation} \label{bias}
\PPP(|\theta(T_0)|\in A)=\int\limits_{0}^{\infty} \PPP(T\in
A~|~l(T)=z)zf_l(z)\ dz,
\end{equation}
where $T$ has under $\PPP$ distribution $m_1^r$.
\end{lemma}

\begin{proof}
Let $F_l$ be the distribution function of $l=l(T)$, and for
$t\in\D{R}$ let $T(t)$ be the 1-slope around $t$. That is, the
slope whose domain contains $t$. Then $\PPP(|\theta(T_0)|\in
A)=\PPP(|\theta(T(t))|\in A)$ for all $t>0$ because $\gq(T_0)$ is
the same as the image under $\gq$ of the slope around $t$ for
$(w_{s-t}-w_{-t}: s\in \D{R})$, and the latter process is again a
standard two sided Brownian motion. Now let $(Y_n)_{n\ge 0}$ be an
independent sequence of slopes with
$Y_n\overset{law}{=}(-1)^{n+1} T$. Glue them sequentially to get
a function $f$ in $C([0,+\infty))$ with $f(0)=0$, and denote by
$\tilde T(t)$ the slope around $t$, for $t>0$.

If we take $\gs$ as defined just before Fact \ref{NPLaplace} with
$x=1$ in the definition of $\tau$ there, then it holds
\begin{equation*}
\PPP(|\theta(T_0)|\in A)=\PPP(|\theta(T(t))|\in
A)=\PPP(|\theta(T(t))|\in A, \gs<t)+\PPP(|\theta(T(t))|\in A,
\gs>t)
\end{equation*}
and
\begin{equation*}
 \PPP(|\gq(T(t))|\in A, \gs <t)=\int_0^t\PPP(|\gq(\tilde
T(t-y))|\in A)f_\gs(y)\ dy,
\end{equation*}
where $f_\sigma$ is the density of $\sigma$. We will take $t\to
+\infty$ and finish with the proof by showing that the limit
$\lim_{t\to+\infty}P(|\gq(\tilde T(t))|\in A)$ exists and
$$\lim_{t\rightarrow+\infty}\PPP(|\gq(\tilde T(t))|\in A)=\int\limits_{0}^{\infty} \PPP(T\in A~|~l(T)=z)zf_l(z)\ dz.$$
To see this, define $g(t):= \PPP(|\gq(\tilde T(t))| \in A)$ for
$t\ge0$. Then
\begin{multline*}
g(t)=\PPP(T\in A,~l(T)> t)+\int\limits_{0}^{t} \PPP(|\gq(\tilde
T(t-s))|\in A)\ dF_{l}(s)\\=\PPP(T\in A,~l(T)> t)+\int
\limits_{0}^{t} g(t-s)\ dF_{l}(s).
\end{multline*}
The distribution of $l$ is nonarithmetic with mean value 1.  By
the renewal theorem (\cite{DU}, Chapter 3, statement (4.9)), it follows that the $\lim_{t\to+\infty} g(t)$
exists and
\begin{eqnarray*}
\lim_{t\to+\infty}g(t)=\int\limits_{0}^{\infty} \PPP(T\in A,~l(T)>
s)\ ds= \int\limits_0^\infty\int\limits_s^\infty \PPP(T\in
A~|~l(T)=z)f_l(z)\ dzds\\=\int\limits_0^\infty\int\limits_0^z
\PPP(T\in A~|~l(T)=z)f_l(z)\ dsdz=\int\limits_{0}^{\infty}
\PPP(T\in A~|~l(T)=z)zf_l(z)\ dz.
\end{eqnarray*}

\end{proof}

Now we apply Lemma \ref{cslope} to obtain the distribution of the
length and height of the central 1-slope.

\medskip

$\bullet$ Regarding the length of $T_0$, observe that for $x>0$ the set $A:=\{T\in\C{S}:l(T)<x\}$ is open and
$$
\PPP(l(T_0)<x)=\int\limits_{0}^{\infty}
\PPP(l(T)<x~|~l(T)=z)zf_l(z)\ dz= \int\limits_0^x z~f_l(z)\ dz.
$$
So that $l(T_0)$ has the density
\begin{equation}
f_{l(T_0)}(x)= xf_l(x) \qquad \text{ for }x\ge0,
\end{equation}
which is the size-biased sampling formula from renewal theory.

\medskip

$\bullet$ Regarding the excess height of $T_0$, observe that for  $x>0$ the set $A:=\{T\in\C{S}:\eta(T)<x\}$ is open
and
$$
\PPP(\eta(T_0)<x)=\int\limits_{0}^{\infty}
\PPP(\eta(T)<x~|~l(T)=z)zf_l(z)\ dz=\int\limits_{0}^{\infty}
\PPP(\eta(T)<x,\ l(T)=z)z\ dz.
$$
Differentiating with respect to $x$ we get the density of
$\eta(T_0)$ as
$$\int\limits_{0}^{\infty}
\PPP(\eta(T)=x,\ l(T)=z)z\ dz =\EEE\left(l(T)~|~\eta(T)=x\right) e^{-x}.
$$
From (iii), (iv) of
Fact \ref{NPLaplace}
$$
\EEE(e^{-tl(T)}~|~\eta(T)=x~)=\psi(t)~e^{-x\phi(t)}.
$$
After some calculations, $\frac{\partial}{\partial
t}\psi(t)~e^{-x\phi(t)}|_{t=0}=-(2x+1)/3$.
\nid The derivative can move inside the expectation on the left hand side of the above relation,
giving $-\EEE(l(T)~|~\eta(T)=x)$, due to the
monotone convergence theorem; which applies because $l(T)>0$ and the function
$(x\mapsto(1-e^{-ax})/x)$ is nonnegative and decreasing in $\D{R}$
for any $a>0$. Therefore the density of $\eta(T_0)$ is
\begin{equation} \label{biasheight}
f_{\eta(T_0)}(x)=\frac{(2x+1)e^{-x}}{3} \qquad \text{ for } x>0.
\end{equation}

The last bit of information needed to achieve our goal, stated at
the first sentence of this section, is the direction of the central 1-slope (i.e., upward or downward) and its
location with respect to zero. By symmetry, $T_0$ is an upward
slope with probability 1/2, and from exercise 3.4.7 of \cite{DU} it
follows easily that given the length $l$ of the slope
$(w-w(x_0))|[x_0,~x_1]$ around zero, the distance of zero from
$x_0$ is uniformly distributed in $[0, l]$.

\section{Proof of the Theorem and the Corollaries}\label{proofs}
For $x>0$ and $w\in \C{W}_1$ with set of $x$-extrema
$R_x(w)=\{x_k:k\in \D{Z}\}$  define
\[
A_x(w):=\Big\{\Big(w-w(x_k)\Big)\big|[x_k,
x_{k+1}]:k\in\D{Z}\Big\}.
\]
We refer to the parameter $x$ as time since we are going to study
the evolution of $A_x(w)$ as $x$ increases. $A_x(w)$ is the set of
slopes at time $x$. Roughly, as $x$ increases, the slopes that have height
smaller than $x$ are absorbed into greater ones.

For any Lebesgue measurable set $S\subset[0,+\infty)$, $x\ge1$,
and $k\in\mathbb{N}$ we define
\[
\C{U}(x, S, k)= \PPP\left( \begin{array}{c} \text{In~~} A_{x}(w)
\text{~~the central slope has \eh ~} y \in S
\\ \text{~~and $b_\cdot(w)$ has changed sign $k$ times in $[1,
x]$~~}\end{array}\right).
\]
It is important to note that the values of $b$ up to time $x$ are
``encoded'' in the central slope of $A_x(w)$. So the number of sign changes of $b_\cdot(w)$ in $[1,x]$
can be inferred by that slope.

For $x, k$ fixed, $\C{U}$ is a measure that satisfies
$$\C{U}(x,S,k) \le \PPP\big( \text{In~~} A_{x}(w) \text{~~the central slope has \eh~} y \in S\big),$$
and the right-hand side, considered as a function of $S$, is a measure absolutely continuous with
respect to the Lebesgue measure with density
$(1/x)f_{\eta(T_0)}(\cdot/x)$, where $f_{\eta_(T_0)}$ is given in
(\ref{biasheight}). Therefore the measure on the left-hand side
has a density also (an element of $L^1([0,+\infty))$),
call it $u(x, y, k)$, and for Lebesgue almost all $y$ it holds
\begin{equation} \label{bound}
\sum\limits_{k=0}^{+\infty} u(x, y, k) =
\frac{1}{x}f_{\eta(T_0)}(\frac{y}{x}).
\end{equation}
Define $U:[1,+\infty)\times[0,
+\infty)\times\D{N}\rightarrow[0,+\infty)$ with
$U(x,y,k):=\C{U}(x,[y,+\infty),k).$ \\ $U$ is continuous as is proved in Lemma \ref{cont} in section \ref{lemmata}. We plan to
establish a PDE for $U$. To do this, we look at $A_x(w)$ and try to
predict how $A_{x+\gep}(w)$ should look like.

\begin{figure}
\begin{center}
\resizebox{8cm}{!}{\input{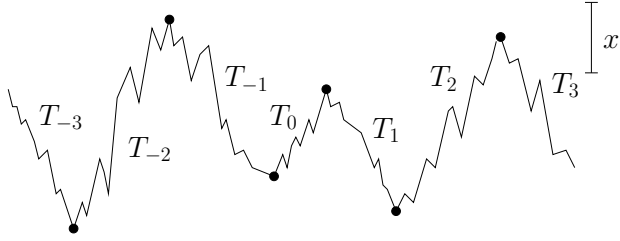}} \caption{\label{P3}
The decomposition of a piece of the Brownian path in $x$-slopes.
The dots mark the points of $x$-extrema. The length $x$ is shown
on the side.}
\end{center}
\end{figure}

In the transition from $A_x(w)$ to $A_{x+\gep}(w)$, a slope around
a point remains the same if the \eh of this and the two neighboring
slopes are greater than $\gep$. In case one slope has \eh in
$[0,\gep)$, it does not appear in $A_{x+\gep}(w)$. For example,
in Figure \ref{P3} the slope $T_0$ has a height say $x+h_0$ with
$h_0\in[0,\gep)$, and the two neighboring slopes $T_{-1}, T_1$,
have height greater than $x+\gep$. Assume that $T_{-1}$, $T_1$,
have heights $x+v_1$,  $x+v_2$ respectively. In $A_{x+\gep}(w)$ we
know that the slopes $T_{-1}, T_0, T_1$  will merge to constitute
a new slope with height $x+v_1+x+v_2-x-h_0=x+v_1+v_2-h_0$; i.e.,
with \eh $v_1+v_2-h_0-\gep$.  The slopes $T_{-2},T_2$ can stay as
they are in the transition from $A_x(w)$ to $A_{x+\gep}$ or they
can be extended if some of $T_{-3}, T_3$ has \eh in $[0,\gep)$. In
any case, they don't interfere with $T_{-1}, T_0,T_1$. This simple
observation, combined with the renewal structure of the sets
$A_x(w)$, is the basis for the next lemma, which is the first step
towards establishing a PDE that $U$ solves. We denote by
$\partial_y U$ the $y$ derivative of $U(x,y,k)$, and recall that
$p_x(v)$ was defined in \eqref{expona} as the density of an
exponential with mean $x$.

\begin{lemma} \label{step1}
For $x\ge1,~ y\ge0,~ \gep>0, k\in\D{N}$
 we have
\begin{multline} \label{preODE}
U(x+\gep,y,k)=U(x,y+\gep,k)e^{-2\gep/x}+\frac{2\gep}{x}\int_0^{+\infty}p_x(v)~U\big(x,(y-v)^+,k\big)~dv-\\
1_{\{k\ge1\}}\gep~\partial_yU(x,0,k-1)(y/x+1)e^{-y/x}+o(\gep),
\end{multline}
where for $k\ge1$ we assume that $U(x,y,k-1)$ is differentiable in
$y$ with continuous derivative. The term $o(\gep)$ depends on
$x,y,k$.
\end{lemma}

\begin{proof}
 The left-hand side of the equation is the probability of an event
referring to $A_{x+\gep}(w)$, and we express it in terms of
probabilities referring to $A_x(w)$. In $A_x(w)$ we focus our
attention on the seven $x$-slopes closest to zero. Denote them by
$T_i,\  i=-3,\ldots,3$ in the order they appear in the path of $w$
from left to right, with $T_0$ being the central one. The slopes
$|\gq(T_i)|$, $i=-3,\ldots,3$ are independent having for $i\neq0$
law $m_x^r$, and $|\gq(T_0)|$ having law $m_x^c$. The probability
of the event that at least two of them have \eh in $[0,\gep)$ is
bounded by $21\gep^2/x^2$, and this is accounted for in the
$o(\gep)$ term in \eqref{preODE}. In the complement of this event,
the event whose probability appears in the left-hand side of
\eqref{preODE} happens if and only if in $A_x(w)$ one of the
following three holds (considering what slope, if any, among the seven has \eh in $[0,\gep)$).
\begin{itemize}
\item At most one of $T_{-3}, T_{-2}, T_2, T_3$ has \eh in
$[0,\gep)$, the excess height of each of $T_{-1}, T_{0}, T_1$ is at least $\gep$,
 and $b$ has changed sign $k$ times in $[1,x]$. In this
case the slope around zero is the same for both
$A_{x}(w),A_{x+\gep}(w)$.
 \item Exactly one of $T_{-1}, T_1$ has \eh in
$[0,\gep)$, $T_0$ has \eh at least $\gep$, and $b$ has changed sign $k$ times in $[1, x]$.
In this case $b_{x+\gep}$ has the same sign as $b_x$. For example, assume that $T_1$ has \eh in $[0,\gep)$.
If $b_x>0$ then $b_{x+\gep}>b_x$, while if $b_x<0$ then $b_{x+\gep}=b_x$ and simply the central slope
in $A_{x+\gep}(w)$ results from merging $T_0, T_1, T_2$.
\item
$k\ge1$, $T_0$ has \eh in $[0, \gep)$, and $b$ has changed sign
$k-1$ times in $[1,x]$. In this case $b_{x+\gep}b_x<0$ because in $A_{x+\gep}(w)$ the central slope
results from merging $T_{-1}, T_0, T_1$ and it has different direction than $T_0$ (i.e., if for example $T_0$
is an upward slope, then the central slope in $A_{x+\gep}(w)$ is a downward slope).
\end{itemize}

\nid The first case has probability
$$U(x,y+\gep,k)\,\PPP(\eta_{-1}\in[\gep,+\infty))\,\PPP(\eta_1\in[\gep,+\infty)) = U(x,y+\gep,k)\exp(-2\gep/x)$$
with $\eta_{-1}$ (resp. $\eta_1$) denoting the \eh of $T_{-1}$
(resp. $ T_1$). This follows by the independence mentioned above
and by the fact that $\eta_{-1}$ and $\eta_1$ have exponential
distribution with mean $x$, and it expresses the demand that both
of the $x$-slopes neighboring the central $x$-slope have \eh
greater than $\gep$. In this case, the excess heights of $T_{-3}, T_{-2}, T_2, T_3$ do not matter.
They can't influence the central slope in $A_{x+\gep}(w)$.

\nid The second case has probability
$$2~\int_0^\gep\int_\gep^{+\infty}p_x(v_1)\,p_x(v_2)\,U\big(x,\gep\vee(y+\gep-v_2+v_1),k\big)\,dv_2dv_1 $$
because, say in the event that $T_1$ has \eh in $[0,\gep)$, the
central slope in $A_{x+\gep}(w)$ comes from merging $T_0,T_1,T_2$
of $A_x(w)$. Assume that they have excess heights $u,v_1,v_2$
respectively. Then the central slope in $A_{x+\gep}(w)$ will have
\eh $u-v_1+v_2-\gep$, and the requirement that this is greater
than $y$ translates to $u$ being greater than
$(y+v_1-v_2+\gep)^+$. And of course, by assumption, $T_0$ and
$T_2$ have \eh greater than $\gep$.

$U$ is continuous as it is proved in Lemma \ref{cont}. Thus,
dividing with $\gep$ the previous integral and taking
$\gep\rightarrow0$ we get as limit
$$\frac{2}{x}\int_0^{+\infty} p_x(v_2)\,U\big(x,(y-v_2)^+,k\big)\,dv_2.$$
From this procedure we pick up another $o(\gep)$ term.

\nid The last case has probability
$$-\int_0^\gep \int_\gep^{+\infty}\int_\gep^{+\infty}\partial_y
U(x,z,k-1)\,p_x(v_1)\,p_x(v_2)\,1_{\{v_1+v_2\ge
y+z+\gep\}}\,dv_1dv_2dz.
$$
By assumption, $-\partial_yU(x,y,k-1)$ exists and it is the density
of the measure $\C{U}(x,\cdot,k)$. The dummy variables $v_1,v_2,z$
stand for the excess heights of $T_{-1},T_1,T_0$ respectively, and
in this case the central slope in $A_{x+\gep}(w)$ has \eh
$v_1+v_2-z-\gep$, giving the restriction $v_1+v_2-z-\gep\ge y$.
Since $\partial_yU(x,\cdot,k-1)$ is continuous, dividing by $\gep$
and taking $\gep\to0$ we get
$$
-\partial_yU(x,0,k-1)\int_0^{+\infty}\int_0^{+\infty}1_{\{v_1+v_2\ge
y\}}\,p_x(v_1)\,p_x(v_2)\,dv_1dv_2.
$$
And again we pick up an $o(\gep)$ term. The double integral equals
$(y/x+1)\exp(-y/x)$.
\end{proof}

\medskip

Before getting to the actual proof of our theorem, we give a
non-rigorous short derivation to illustrate its main steps. The
main problem is that we don't know if $U$ is differentiable in the
$x,y$ variables for every $k\in\D{N}$. Assume for the moment that
it is. For $f,g\in L^1\big([0,+\infty)\big)$, as usual, we define
$f\ast g\in C\big([0,+\infty)\big)$ by $(f\ast g)(x):=\int_0^x
f(x-y)g(y)\,dy$ for all $x\in[0,+\infty)$.

The above lemma would give for $U$ the PDE

\begin{multline}\label{firstPDES}
(\partial_x-\partial_y)\,U(x,y,k)=-\frac{2}{x}~U(x,y,k)+
\frac{2}{x^2}~\big(U(x,\cdot,k)\ast
e^{-\frac{\cdot}{x}}\big)(y)\\+\frac{2}{x}e^{-\frac{y}{x}}U(x,0,k)-1_{\{k\ge1\}}
\partial_yU(x,0,k-1)(y/x+1)\exp(-y/x).
\end{multline}
Let $f(x,y,k)=U(x, yx,k)$; i.e., $U(x,y,k)= f(x, y/x,k).$
\\Then $f$ should satisfy
\begin{multline}\label{finalPDES}
(x\,\partial_x-(1+y)\,\partial_y+2)f(x,y,k)=2\,\big(f(x,\cdot,k)\ast e^{-\cdot}\big)(y)+2\,e^{-y}f(x,0,k)\\
-1_{\{k\ge1\}}\,(y+1)\,e^{-y}\,\partial_y f(x,0,k-1),
\end{multline}
while the conditions at $x=1$ are
\begin{align}
f(1,y,0)&=(2y/3+1)\,e^{-y}& \text{for~~} y\ge0, \label{PDEconditions0}\\
f(1,y,k)&=0 & \text{for~~} y\ge0, k\ge 1. \label{PDEconditionsk}
\end{align}
The first equation comes from \eqref{biasheight}, the second is
clear.

For $z\in D:=\{z\in\D{C}:|z|<1\}$, the generating function
$M(x,y,z):=\sum\limits_{k=0}^{\infty}f(x,y,k)\, z^k$ is well
defined. Assuming that $M$ is differentiable with respect to $x,y$
and its $x,y$ derivatives are obtained with term by term
differentiation, we see that $M$ satisfies the PDE problem
\begin{align}\label{MPDE}
(x\,\partial_x-(1+y)\,\partial_y+2)\,M(x,y,z)&=2\,\big(M(x,\cdot,z)\ast e^{-\cdot}\big)(y)\notag\\
+2\,e^{-y} M(x,0,z)-&(y+1)\,e^{-y}\,z\partial_y M(x,0,z) \text{\quad in } (1,\infty)\times(0,\infty),\\
M(1,y,z)&=(2y/3+1)\,e^{-y} \text{\qquad  for \quad} y\ge0.
\label{Minitial}
\end{align}
We try for a solution of the form
$$M(x,y,z)=[a(x,z)+b(x,z)\,y]\,e^{-y}.$$
Substituting this into \eqref{MPDE} we see that $e^{-y}$ factors
out in both sides, and after cancellation, we arrive in an
equality of two first degree polynomials in $y$ with coefficients
depending on $x,z$. Equating the coefficients in equal powers of
$y$ in the two sides of the equation, we arrive at the following
system of ODEs for $a,b$.
\begin{align*}
x\,\partial_x\, a(x,z)+(z-1)\,\left(b(x,z)-a(x,z)\right)&=0,\\
x\,\partial_x \,b(x,z)+(2+z)\,b(x,z)-(z+1)\,a(x,z)&=0.
\end{align*}
The initial condition \eqref{Minitial} for $M$, expressed in terms
of $a, b$, becomes $a(1,z)=1,\, b(1,z)=2/3$ for all $z\in D$. We
easily see that the only solution of the system satisfying these conditions is
\label{alfa}
\begin{align*}
a(x,z)&=c_1(z)\,x^{\gl_1(z)}+c_2(z)\,x^{\gl_2(z)},   \\
b(x,z)&=c_1(z)\left(1+\frac{\gl_1(z)}{1-z}\right)x^{\gl_1(z)}+c_2(z)\left(1+\frac{\gl_2(z)}{1-z}\right)x^{\gl_2(z)},
\end{align*}
where
$$\gl_{1}(z)=\frac{-3+\sqrt{5+4z}}{2},\quad
\gl_{2}(z)=\frac{-3-\sqrt{5+4z}}{2},$$ and
$$c_{1}(z)=\big((z-1)/3-\gl_2(z)\big)/(\gl_1(z)-\gl_2(z)), \quad
 c_{2}(z)=\big(-(z-1)/3+\gl_1(z)\big)/(\gl_1(z)-\gl_2(z)).$$
Now $E(z^{k(x)})=M(x,0,z)=a(x,z)$, which is what we want.

\medskip

\textbf{Proof of the Theorem:} The proof is done by taking
the steps of the above ``proof'' in reverse order. This time all
the steps can be justified. We will need three lemmata. Two of them are non-trivial, and
their proof is given in section \ref{lemmata}.

\begin{lemma}\label{powerseries}
The solution $M$ of \eqref{MPDE}, \eqref{Minitial} obtained above
is analytic as a function of $z$ in $D$. The coefficients of its
development as a power series around zero are differentiable with
respect to $x$ and $y$ in $(1,+\infty)\times(0,+\infty)$, and its
$x,y$-derivatives are continuous and can be found with term by
term differentiation.
\end{lemma}

\nid Using this lemma, we write $M$ as
\[
M(x,y,z):=\sum\limits_{k=0}^{\infty}g(x,y,k)\,z^k.
\]
Differentiating $M$ term by term and equating the coefficients of
equal powers of $z$ in the two sides of \eqref{MPDE}, we see that
the sequence of functions $\big(g(\cdot,\cdot,k)\big)_{k\ge 0}$
satisfies the PDEs \eqref{finalPDES} with conditions at $x=1$
given by \eqref{PDEconditions0}, \eqref{PDEconditionsk}.

For $k\in\D{N}$ define $\tilde
g(\cdot,\cdot,k):[1,+\infty)\times[0,+\infty)\to\D{R}$ by $\tilde
g(x,y,k)=g(x,y/x,k)$ for $(x,y,k)\in[1,+\infty)\times[0,+\infty)$.
The sequence of functions $\big(\tilde g(\cdot,\cdot,k)\big)_{k\ge
0}$ satisfies the PDEs \eqref{firstPDES} with conditions at $x=1$
given by \eqref{PDEconditions0}, \eqref{PDEconditionsk}.

The proof is finished by showing that the sequence
$\big(U(\cdot,\cdot,k)\big)_{k\in\D{N}}$ satisfies a weak form of
these PDEs, and then a uniqueness result will identify $U$ as
$\tilde g$.

For fixed $c>1$, define $g_{c,k}(x)=U(x,c-x,k)$  for $x\in[1, c].$
We state as a lemma an equation that $g_{c,k}$ satisfies. The
proof is straightforward from Lemma \ref{step1}.

\begin{lemma} \label{ODES}
The function $g_{c,k}$ is differentiable in $(1, c)$ and satisfies
\begin{multline} \label{ODE}
g_{c,k}^\prime(x)= -\frac{2}{x}~g_{c,k}(x)+
\frac{2}{x^2}~\big(U(x,\cdot,k)\ast
e^{-\frac{\cdot}{x}}\big)(c-x)\\
+\frac{2}{x}\,e^{-\frac{c-x}{x}}\,U(x,0,k)-1_{\{k\ge1\}}\,\partial_yU(x,0,k-1)\,(c/x)\,e^{-(c-x)/x},
\end{multline}
where for $k\ge1$ we assume that $U(x,y,k-1)$ is differentiable in
$y$ with continuous derivative.
\end{lemma}
\nid And the promised uniqueness result is the following.
\begin{lemma}[Uniqueness]  \label{uniq}
Let $f:[1,+\infty)\times[0, +\infty)\rightarrow[0,+\infty)$ be a
continuous function such that for every $c>1$ the function
$g_c:[1,c]\rightarrow\mathbb{R}$ with $g_c(x):=f(x, c-x)$ is
continuous on $[1,c]$ and differentiable on $(1, c)$, with
\begin{equation} \label{PDE}
g_c'(x)= -\frac{2}{x}~g_c(x)+ \frac{2}{x^2}~(f(x,\cdot)\ast
e^{-\frac{\cdot}{x}}\big)(c-x)+\frac{2}{x}e^{-\frac{c-x}{x}}f(x,0)
\end{equation}
and $g_c(1)=0$. Then $f\equiv0$.
\end{lemma}
\nid Now using induction we show that
\[
U(\cdot,\cdot,k)=\tilde g(\cdot,\cdot,k) \qquad \forall k
\in\D{N}.
\]
The function $U(\cdot,\cdot,0)-\tilde g(\cdot,\cdot,0)$ satisfies
the assumptions of Lemma \ref{uniq} because of the PDE problem
that $\tilde g(\cdot,\cdot,0)$ solves and Lemma \ref{ODES}. For
$k\ge1$, assuming the statement true for $k-1$, the same argument
works for $U(\cdot,\cdot,k)-\tilde g(\cdot,\cdot,k)$, where now
the assumption on $U(\cdot,\cdot,k-1)$ required by Lemma
\ref{ODES} is provided by the inductive hypothesis.

Therefore, $\sum_{k=0}^{+\infty}U(x,xy,k)\,z^k=M(x,y,z)$ for
$(x,y,z)\in[1,+\infty)\times[0,+\infty)\times D$ and
$\EEE(z^{k(x)})=M(x,0,z)=a(x,z)$, proving \eqref{genfun} for $z\in
D$.  \hfill \qedsymbol

\medskip

\textbf{Proof of Corollary \ref{Markov}:}  Our theorem gives
that $\PPP(b$ doesn't change sign in $[1,x])=a(x,0)$, where the
function $a$ is defined on page \pageref{alfa}, and by scaling
\[\PPP(\text{$b$ doesn't change sign in
}[x, y])=a(y/x,0) \text{\qquad for $0<x<y$.}
\]
Consequently, the density of the last point before $y$ that we
have sign change is  $-y x^{-2}\partial_x a(y/x, 0)$, where we use
$\partial_x$ here and below to denote derivative with respect to
the first argument. Differentiating with respect to $y$, we get
the density of the event that $x, y$ are consecutive times of sign
change as
\[x^{-2}\partial_x a(y/x, 0)+y x^{-3}\partial_{xx}a(y/x,0),\]
which, after using the expression for $a(x,0)$, becomes
\[\frac{1}{3(\gl_1-\gl_2)xy}\left((y/x)^{\gl_1}-(y/x)^{\gl_2}\right).
\]
The event that $b$ changes sign at $x$ translates to the central
$x$-slope having excess height 0. This has density $1/(3x)$ due to
the scaling property of the $x$-slopes and relation
\eqref{biasheight}, which refers to the central 1-slope. Thus, the
density of the time $Y$ where the next sign change after $x$
happens, given that there was a sign change at $x$, is
\begin{equation}\label{transition}
h(y)=\begin{cases}\frac{1}{(\gl_1-\gl_2)y}\big((y/x)^{\gl_1}-(y/x)^{\gl_2}\big)
& \text{if } y\ge x,\\
0 & \text{otherwise}.
\end{cases}
\end{equation}
The quotient $Y/x$, given that there was a sign change at $x$, has
the density of the $r_i$'s given in \eqref{rdensity}.

 \hfill \qedsymbol

\textbf{Proof of Corollary \ref{strenghening}:} Observe that the
right-hand side of \eqref{genfun} is a function analytic in
$\D{C}\smallsetminus(-\infty,-5/4].$ As for the left-hand side, we
have the following lemma.
\begin{lemma} \label{entire}
For any $x>1$, the power series
$\sum_{n=0}^{+\infty}\PPP(k(x)=n)~z^n$ defines an entire function.
\end{lemma}
\nid And from a basic property of analytic functions, it follows
that the quantities $M(x,0,z)$, $a(x,z)$ agree for all
$z\in\D{C}\smallsetminus(-\infty,-5/4]$.  \hfill \qedsymbol

\begin{remark}
Of course, Lemma \ref{entire} implies that the right hand
side of \eqref{genfun} can be extended to an entire function.
However, the way \eqref{genfun} is written doesn't allow as to claim that it
holds for all $z\in\D{C}$ because the function $(z\to\sqrt{4+4z})$ doesn't have an entire extension.
\end{remark}

\textbf{Proof of Corollary \ref{LDP1}:} We apply the Gartner-Ellis
Theorem  (Theorem 2.3.6 in \cite{DZ}). The moment generating
function of $k(e^t)/t$ is given for any $\lambda\in\D{R}$ by
$\Lambda_{t}(\lambda):=\log \EEE(\exp\{\lambda k(e^t)/t\}) $. So
$t^{-1}\Lambda_t(t\lambda)=t^{-1}\log \EEE(e^{\lambda k(e^t)})=
t^{-1}\log M(e^t,0, e^\lambda)$, and using Corollary
\ref{strenghening} we see that
$$\Lambda(\lambda):=
\lim\limits_{t\rightarrow\infty}\frac{1}{t}\Lambda_{t}(t\lambda)=\lambda_1(e^\lambda)=\frac{-3+\sqrt{5+4e^\lambda}}{2}.
$$
The Fenchel-Legendre transform $\gL^*$ of $\gL$, defined by
$\gL^*(x)=\sup_{\gl\in\D{R}}\{\gl x-\gL(\gl)\}$ for all
$x\in\D{R}$, is found to be the function $I$ defined in the
statement of the proposition. Also
$D_\Lambda:=\{\lambda\in\mathbb{R}:\Lambda(\lambda)<\infty\}=\mathbb{R},$
and $\Lambda$ is strictly convex and differentiable in
$D_\Lambda$. The result follows from the Gartner-Ellis theorem.
\hfill  \qedsymbol

\section{Proofs of the lemmata}\label{lemmata}

\begin{lemma} \label{cont}
For any $x\ge1$, $y\ge0,$ $k\in\mathbb{Z}$, and $\gep>0,$ we have
\begin{align*}
|U(x+\gep, y, k)-U(x, y+\gep, k)|&\le3\gep/x,\\
|U(x, y, k)-U(x, y+\gep, k)|&\le\gep/x.
\end{align*}
In particular, $U$ is continuous.
\end{lemma}

\begin{proof}
Call $A$ and $B$ the two events whose probabilities are $U(x+\gep,
y, k)\text{~~and~~}U(x, y+\gep, k)$ respectively. Then $A\triangle
B\subset$ [\ In $A_x(w)$ at least one of the three slopes
neighboring zero has excess height $<\gep$\ ]. Denote by $T_0,
T_1$ the central slope and the slope to the right to it in
$A_x(w)$. Then \begin{multline*} \PPP(A\triangle B)\le 2\PPP(\ T_1
\text{~~has excess height~~} <\gep \ )+\PPP(\ T_0 \text{~~has
excess height~~} <\gep \ )
\\=2\int_0^{\gep/x}e^{-z}\,dz+\int_0^{\gep/x}(2z/3+1/3)e^{-z}\,dz<3\gep/x.
\end{multline*}
The other inequality follows because the density of the measure
$\mathcal{U}(x,S,k)$ is bounded by $1/x$ (see (\ref{bound})).
\end{proof}

\textbf{Proof of Lemma \ref{powerseries}:}
 Define
$K:\D{C}\times\D{C}\times D\to \D{C}$ with $K(z_1, z_2,
z):=M(e^{z_1}, z_2, z)$ for $(z_1, z_2, z) \in
\D{C}\times\D{C}\times D$. Clearly, $K$ is a holomorphic function
(choose an analytic branch of the square root function defined on
$\D{C}\setminus(-\infty,0)$; the number $5+4z$ is there for $z\in
D$) and has a power series development centered at zero that
converges in $\D{C}\times\D{C}\times D$ (see, e.g., \cite{KR},
Proposition 2.3.16). The claims of the lemma follow by the
relation $M(x, y, z)=K(\log x, y, z)$ and standard properties of
power series.
 \hfill \qedsymbol

\textbf{Proof of Lemma \ref{uniq}:} For $c>1$ and $x\in[1,c]$,
define $N(c,x):=\sup\{|f(z, c-z)|: z\in[1,x]\}$. $N$ is well
defined because $f$ is bounded on compact sets.
\\ From (\ref{PDE}), for $x\in(1,c)$ one has
$$|g_c'(x)|\le2~N(c,x)+4~\sup_{1\le d\le c}N(d,x),$$
and since $g_c(1)=0$, we get after integrating
$$ |f(x, c-x)|\le 6\int_1^x \sup_{1\le d\le c}N(d,t)dt,$$
which implies
$$N(c,x)\le 6\int_1^x \sup_{1\le d\le c}N(d,t)dt$$
and
$$ \sup_{1\le d\le c}N(d,x)\le 6\int_1^x \sup_{1\le d\le c}N(d,t)dt.$$
The function $A(x)=\sup_{1\le d\le c}N(d,x)$ ($x\in [1,c]$) is
continuous (because $f$ is) and has $A(1)=0$.  An application of
Gronwall's lemma to $A$ gives $N(d,x)=0$ for $x\in[1,c],\
d\in[1,c]$; i.e., $f(x,y)=0$ for all $(x,y)\in[1,+\infty)\times[0,
+\infty)$. \hfill \qedsymbol

\medskip

\textbf{Proof of Lemma \ref{entire}:} Let $X_k$ be as in Corollary
\ref{Markov}. Then
\begin{multline*} \PPP(k(x)=n+1)\le\PPP(k(x)\ge n+1)=\PPP(X_{n+1}\le
x) \\
\le \PPP(r_1r_2\cdots r_n\le x)=\PPP(\log r_1+\log r_2+ \cdots
+\log r_n\le \log x).
\end{multline*}
For $i\ge 1$ set $Y_i=\log r_i$, $S_i=Y_1+Y_2+\cdots+Y_i$, and let
$m_i$ be the distribution measure of $S_i/i$. By Cramer's theorem
(\cite{DZ} Theorem 2.2.3), the sequence $(m_i)_{i\ge1}$ satisfies
a large deviation principle with rate function
$I(x)=\sup_{\gl\in\D{R}}\{\gl x-\log \EEE(e^{\gl \log r_1})\}$ for
$x\in\D{R}$.

\nid Clearly, for any $\gep>0$ we have
$$\varlimsup_{n\to\infty} \frac{1}{n}\log \PPP(\frac{S_n}{n}<\frac{\log
x}{n})<-I(\gep).$$
And $\lim_{\gep\to 0}I(\gep)=+\infty$ because
$I$ is lower semicontinuous and $I(0)=+\infty$. To see the last
point, observe that $$I(0)=\sup_{\gl\in\D{R}}\{-\log \EEE(e^{\gl
\log r_1})\}\ge\varlimsup_{\gl\to -\infty}\{-\log \EEE(e^{\gl \log
r_1})\}=+\infty$$ because $\lim_{\gl\to -\infty}e^{\gl \log r_1}=0$
with probability one and the bounded convergence theorem applies.
Thus $\lim_{n\to\infty}\frac{1}{n}\log \PPP(k(x)=n+1)=-\infty$,
proving that the radius of convergence for the power series is
infinite. \hfill \qed

\begin{lemma}\label{discrete}
For $\PPP$ and $\C{W}_1$ as defined in the introduction, it holds that $\PPP(\C{W}_1)=1$.
\end{lemma}
\begin{proof}
First we prove that for fixed $x>0$, the set $$C_x:= \{w\in C(\D{R}): \text{ $R_x(w)$ has the properties appearing in the definition of $\C{W}_1$}\}$$ has $\PPP(C_x)=1$. Observe
that for $z$ a point of $x$-minimum and
$\ga_z:=\sup\{\ga<z:w(\ga)\ge w(z)+x\}$,
$\gb_z:=\inf\{\gb>z:w(\gb)\ge w(z)+x\}$ it holds that
$\ga_z<z<\gb_z$ because $w$ is continuous at $z$. And there is no other point of $x$-minimum
in $(\ga_z, \gb_z)$ since if $\tilde z$ is such a point, say in $(\ga_z, z)$,
then assuming $\gb_{\tilde z}<z$  we get a contradiction with the definition of $a_z$
while assuming $\gb_{\tilde z}>z$ we get that $w$ takes the same value in two local minima,
which has probability zero. Now assume that there is a strictly
monotone, say increasing, sequence $(z_n)_{n\ge1}$ of $x$-minima converging to
$z_\infty\in\D{R}$. By the above observations we get
$\varlimsup_{y, \tilde y\nearrow z_\infty}(w(y)-w(\tilde y))\ge x$
 contradicting the continuity of $w$ at $z_{\infty}$. Similarly if
$(z_n)_{n\ge1}$ is decreasing. So in a set of $w$'s in $C(\D{R})$ with
probability 1, it holds that the set of $x$-minima of $w$ has no
accumulation point. The same holds for the set of $x$-maxima and
as a result also for $R_x(w)$. Since $\varliminf_{|t|\to\infty}
w_t=-\infty$, $\varlimsup_{|t|\to\infty} w_t=+\infty$,  it follows
that $\PPP(R_x(w) \text{ is unbounded above and below} )=1$. It is clear that
between two consecutive $x$-maxima (resp. minima) there is exactly
one $x$-minimum (resp. $x$-maximum). Consequently, $\PPP(C_x)=1$.

Finally, note that for all $n\in\D{N}\setminus \{0\}$ we have $R_n(w)\subset R_x(w)\subset
R_{1/n}(w)$ for $x\in [1/n, n]$, from which it follows that
$\C{W}_1=\cap_{x\in(0, +\infty)} C_x=\bigcap_{n\in\D{N}\setminus\{0\}}(C_n\cap C_{1/n})$. Thus, $\PPP(\C{W}_1)=1$.
\end{proof}

\bigskip

\textbf{Acknowledgments.} I thank my advisor Amir Dembo for his
help throughout the time this research was conducted and for his
comments on previous drafts of the paper. I also thank Ofer
Zeitouni for his comments on an earlier draft, Daniel Fisher for a
useful discussion we had during a visit of his to Stanford, and
Francis Comets for providing me with a preprint of \cite{CP}.

\end{document}